\documentclass[1p]{elsarticle}
\usepackage{amssymb}
\usepackage{amsfonts}
\usepackage{amsmath}

\usepackage[polish, english]{babel}
\usepackage{polski}
\usepackage{tikz}

\newtheorem{theorem}{Theorem}
\newtheorem{conjecture}[theorem]{Conjecture}
\newtheorem{corollary}[theorem]{Corollary}
\newtheorem{lemma}[theorem]{Lemma}
\newtheorem{observation}[theorem]{Observation}
\newtheorem{problem}[theorem]{Problem}

\newproof{pf}{Proof}

\begin{document}

\title{On asymptotically tight bounds for the open conflict-free chromatic indexes of nearly regular graphs} 
\tnotetext[t1]{This research was supported by the AGH University of Krakow under grant no. 16.16.420.054, funded by the Polish Ministry of Science and Higher Education.}

\author[agh]{Mateusz Kamyczura}
\ead{kamyczuram@gmail.com}

\author[agh]{Jakub Przyby{\l}o}
\ead{jakubprz@agh.edu.pl}

\address[agh]{AGH University of Krakow, al. A. Mickiewicza 30, 30-059 Krakow, Poland}

\begin{abstract}

An edge colouring $c$ of a graph $G$ is called  \emph{conflic--free} if every non--isolated edge of $G$
has a uniquely coloured neighbour in its open edge neighbourhood.
The least number of colours admitting such a colouring is denoted by $\chi'_{\rm OCF}(G)$, or $\chi'_{\rm pOCF}(G)$ if we additionally require $c$ to be proper.

Our main result implies in particular that $\chi'_{\rm OCF}(G) \le \log_2 \Delta + O(\ln\ln\Delta)$
for nearly regular graphs $G$ with maximum degree $\Delta$, which is asymptotically optimal, as witnessed by
the complete graphs.
For proper colourings, we moreover show that $\chi'_{\rm pOCF}(G) \le \Delta + O(\ln \Delta)$ in the same
regime.
These results 
improve existing bounds stemming from related colouring models
and transfer directly to random graphs' setting.

The proofs combine decomposition techniques with probabilistic arguments and
structural properties of edge neighbourhoods.
\end{abstract}

\begin{keyword}
conflict‐free colouring \sep conflict‐free index \sep conflict‐free number  \sep open edge neighbourhood
\end{keyword}

\maketitle

\section{Introduction}

Let $G=(V,E)$ be a finite simple graph with maximum degree $\Delta$ and minimum
degree $\delta$.
By an \emph{edge colouring} of $G$ we mean any mapping
$c : E \to \mathbb{N}$, unless stated otherwise; the colouring is not required to be proper.
For a vertex $v\in V$, we denote by $E_G(v)$ the set of all edges incident with $v$ in $G$.
For an edge $e=uv\in E$, we define its \emph{open neighbourhood} as
\[
E_G(uv) := \bigl(E_G(u)\cup E_G(v)\bigr)\setminus\{uv\}.
\]

We say that an edge $e$ is \emph{satisfied} by a colouring $c$ if among the colours
appearing on the edges in $E_G(uv)$ there exists at least one colour that occurs
exactly once.
If every non-isolated edge of $G$ is satisfied, then $c$ is called an
\emph{open conflict--free edge colouring} of $G$.
The minimum number of colours in such a colouring is the
\emph{open conflict--free chromatic index}, denoted by $\chi'_{\rm OCF}(G)$.
If the colouring is additionally required to be proper, we write
$\chi'_{\rm pOCF}(G)$ to denote the corresponding modified graph invariant.
Clearly,
$\chi'_{\rm OCF}(G)\le \chi'_{\rm pOCF}(G)$,
and, as we shall see, these two parameters 
exhibit fundamentally different behaviour.

The conflict--free paradigm was initially geometrically  motivated -- by frequency assignment problems in cellular
networks and wireless communication, where a base station must serve clients in its vicinity 
without interference~\cite{EvenEtAl,KostochkaEtAl,SmorodinskyPhd,SmorodinskyApplications}.

The related early graph--theoretic research 
focused on \emph{vertex} colourings with a unique colour
in every \emph{closed} neighbourhood.
Let us denote by  $\chi_{\rm CF}(G)$ the corresponding graph invariant.
Foundational results of Glebov, Szab\'o and Tardos~\cite{GlebovEtAl}, as well as of
Bhyravara\-pu, Kalyanasundaram and Mathew~\cite{Hindusi}, established the asymptotic
behaviour of the conflict--free chromatic number in terms of $\Delta$.
In particular, these results imply that there exist graphs with
$\chi_{\rm CF}(G)=\Omega(\ln^2\Delta)$, while, on the other hand,
$\chi_{\rm CF}(G)=O(\ln^2\Delta)$ holds in general.

The variant based on open neighbourhoods turned out to be significantly more
challenging.
For instance, subdivisions of complete graphs show that one may have
$\chi_{\rm OCF}(G)=\Delta+1$ for the corresponding open neighbourhood parameter.
Nevertheless, Pach and Tardos~\cite{PachTardos} proved that
\begin{equation}\label{Eq1.1}
\chi_{\rm OCF}(G)=O(\ln^2\Delta)
\end{equation}
for graphs with minimum degree $\delta=\Omega(\ln\Delta)$.
Furthermore, results of Bhyravarapu, Kalyanasundaram and Mathew~\cite{Hindusi3,Hindusi2}
imply in particular that
\begin{equation}\label{Eq1.2}
\chi_{\rm OCF}(G)=O(\ln^{2+\varepsilon}\Delta)
\end{equation}
for all claw--free graphs~$G$.

In recent years, proper conflict--free colourings and their variants
have been studied intensively for many graph classes, including sparse,
planar, and minor--closed families; see, for example,
\cite{TreesCF,OuterplanarCF,CPS,EUN,Fab,DWC-CHL}.
More recent research have developed this line of research towards the theory of
proper $h$--conflict--free colourings,
introduced by Cho, Choi, Kwon and Park~\cite{hCFBounds} as a natural generalisation
interpolating between conflict--free and distance--two colourings.
A proper vertex colouring $c$ of a graph $G$ is called $h$--conflict--free if the open neighbourhood $N(v)$ of every vertex $v$ hosts $\min\{d(v),h\}$ distinct colours each of which occurs exactly once in $N(v)$ (note that for $h\geq \Delta-1$, such $c$ is a distance--two colouring of $G$).
We denote by $\chi^h_{pOCF}(G)$ the least number of colours admitting such a colouring.
In~\cite{CF-Pirot}, improving earlier results of Liu and Reed~\cite{Liu-Reed} for the case $h=1$,
Chuet, Dai, Ouyang and Pirot obtained asymptotically sharp bounds of the
form
\begin{equation}\label{Eq1.3}
\chi^h_{\rm pcf}(G) \le h\Delta + O(\log\Delta),
\end{equation}
for every fixed $h$,
and conjectured that the bound $h\Delta+1$ should hold for sufficiently large $\Delta$.
A closely related notion of \emph{peaceful colourings}, which further highlights
the connection between conflict--free constraints and distance--based colouring
models, was introduced by Liu and Reed~\cite{Liu-Reed_2}.
In particular, their results strongly support the above conjecture of
Chuet, Dai, Ouyang and Pirot for nearly regular graphs and relatively small values of~$h$,
by showing that
\begin{equation}\label{Eq1.4}
\chi^h_{\rm pCF}(G) \le \Delta + 1,
\end{equation}
whenever $\delta > 8000\Delta/8001$ and $h \le \delta - 8000\Delta/8001$.
This result also implies the conjecture of Caro, Petru\v{s}evski and \v{S}krekovski
from~\cite{CPS} for the corresponding class of graphs.

In contrast to the extensive literature devoted to conflict--free \emph{vertex}
colourings in the open--neighbourhood setting,
the corresponding notion of conflict--free \emph{edge} colourings has received
almost no direct attention.
Most existing results on conflict--free edge colourings concern the closed--neighbourhood
model.
In this framework, D\k{e}bski and Przyby{\l}o~\cite{DebskiPrzybylo} proved an
order--wise optimal upper bound for the corresponding parameter
(for not necessarily proper edge colourings),
namely $\chi'_{\rm CF}(G)=O(\log_2\Delta)$.
Subsequent research focused on refining the leading constant in this bound,
in particular for nearly regular and random graphs; see~\cite{MPK,KPreg,GuoTrees24}.
(Note that in this model, imposing properness would reduce the problem to
standard proper edge colourings.)

The open--neighbourhood edge variant studied here exhibits behaviour that is
substantially different from that of $\chi'_{\rm CF}(G)$ and requires 
development of new ideas.
In this paper, we initiate a systematic study of this model.
Clearly, studying $\chi'_{\rm OCF}(G)$ is equivalent to investigating
conflict--free vertex colourings of the line graph $L(G)$.
Since line graphs are claw--free, the bound~(\ref{Eq1.2}) immediately implies that
for every graph $G$ with maximum degree $\Delta$,
$
\chi'_{\rm OCF}(G)=O(\ln^{2+\varepsilon}\Delta)
$
for any fixed $\varepsilon>0$.
Moreover, by~(\ref{Eq1.1}), this bound improves to
$\chi'_{\rm OCF}(G)=O(\ln^{2}\Delta)$ in the nearly regular regime.
We in particular prove that
$
\chi'_{\rm OCF}(G)\le (1+o(1))\log_2\Delta
$
for nearly regular graphs with maximum degree $\Delta$; see
Theorem~\ref{Main1} in Section~\ref{FinalResultsSection}.
This result is asymptotically optimal, since, for instance, the complete graphs
provide examples with
$\chi'_{\rm OCF}(G)\ge (1-o(1))\log_2\Delta$
(see the discussion in Section~\ref{ConcludingRemarksSection}).
For proper colourings, we additionally show that
$\chi'_{\rm pOCF}(G)\le \Delta+O(\ln\Delta)$ for nearly regular graphs; see
Theorem~\ref{Main2} in Section~\ref{FinalResultsSection}.
By contrast, the general bound obtained from~(\ref{Eq1.3}) yields only
$\chi'_{\rm pOCF}(G)\le 2\Delta+O(\ln\Delta)$,
with a slight improvement in the nearly regular case following from~(\ref{Eq1.4}).

To achieve these results, we develop several new ideas and combine them with
adaptations of existing tools to the open--neighbourhood setting, as well as
e.g. observations originating in research on \emph{majority edge colourings}.
In particular, we employ a modified version of a decomposition lemma introduced
in~\cite{KPreg} (see Section~\ref{MainLemmaSection}), together with a refined
Talagrand--type concentration argument inspired by~\cite{DebskiPrzybylo}
(see Section~\ref{MainResultSection}).
As an application, we derive corresponding results for random graphs $G(n,p)$ in
Section~\ref{FinalResultsSection}, where we summarise and formally state most of our results.
The concluding Section~\ref{ConcludingRemarksSection}
includes
several comments concerning our approach and results, supplemented with suggestions of further research and open problems.

\section{Tools}

Given a graph $G=(V,E)$, subsets $V'\subseteq V$, $E'\subseteq E$ and a vertex $v\in V$, we denote by
$d_{V'}(v)$ the number of edges joining $v$ with $V'$, and we denote by $d_{E'}(v)$ the number of edges in $E'$ incident with $v$. An edge-decomposition of $G$ is in turn a collection of subgraphs $G_1,\ldots,G_k$ of $G$ whose edge sets partition $E$; for definiteness we also assume that $V(G_i)=V$ for every $i\in[k]$. 
Finally, for disjoint subsets $A,B\subseteq V$, let $G[A,B]$ denote the bipartite subgraph of $G$ with the partition sets $A$, $B$ and the edge set consisting in all the edges of $G$ joining $A$ with $B$.

We shall rely on a collection of standard probabilistic tools,
starting with the well-known symmetric variant of the Lov\'asz Local Lemma; 
see e.g.~\cite{LLL}.

\begin{lemma}[Lov\'asz Local Lemma]\label{LLL}
Let $\Omega$ be a finite family of events in a probability space.  
Assume that every event $A\in\Omega$ is mutually independent of all but at most $D$ other events in $\Omega$, 
and that $\Pr(A)\le p$ for each $A\in\Omega$.  
If
\[
e p (D+1) \le 1,
\]
then $\Pr\!\left(\bigcap_{A\in\Omega}\overline{A}\right)>0$.
\end{lemma}

We shall also make frequent use of Chernoff-type concentration bounds.  
The two variants below shall be sufficient for all our estimates; see e.g.~\cite{RandOm,ChernoffBook}.

\begin{lemma}[Chernoff Bound I]\label{Ch}
Let $X=\sum_{i=1}^n X_i$ where the $X_i$ are independent Bernoulli variables with 
$\Pr(X_i=1)=p_i$.  Then for any $0<t\le \mathbb{E}(X)$,
\[
\Pr\!\left(X \le \mathbb{E}(X)-t\right)
   \le \exp\!\left(-\frac{t^2}{2\mathbb{E}(X)}\right).
\]
\end{lemma}

\begin{lemma}[Chernoff Bound II]\label{Ch2}
Under the same assumptions as above, for every $t>0$,
\[
\Pr\!\left(X \ge \mathbb{E}(X)+t\right)
   \le \exp\!\left(-\frac{t^2}{t+2\mathbb{E}(X)}\right),
\]
and in particular,
\[
\Pr\!\left(X \ge \mathbb{E}(X)+t\right)
   \le \exp\!\left(-\frac{t^2}{3\mathbb{E}(X)}\right)
   \quad 
   for \quad 0<t\le\mathbb{E}(X).
\]
\end{lemma}

Note that both inequalities remain valid (with minor rephrasing) whenever only suitable bounds on 
$\mathbb{E}(X)$ are known, rather than its precise value.

We shall also need Talagrand's Inequality -- a concentration tool applicable to
random variables which depend in a controlled manner on a list of independent trials.  
We state below its formulation from~ \cite{MolloyTalagrandRef}, which in convenient for our application.

\begin{theorem}[Talagrand's Inequality]\label{tlg}
Let $X$ be a nonnegative random variable determined by $\ell$ independent trials 
$T_1,\dots,T_\ell$.  
Assume there exist constants $c,k>0$ such that:
\begin{enumerate}
    \item changing the outcome of any single trial alters $X$ by at most $c$, and
    \item for each $s>0$, if $X\geq s$ then there is a set of at most $ks$ trials whose outcomes certify that $X\geq s$.
\end{enumerate}
Then for any $t\ge 0$,
\[
\Pr\!\left(|X-\mathbb{E}(X)| 
        > t + 20c\sqrt{k\mathbb{E}(X)} + 64c^2 k \right)
   \le 4\exp\!\left(-\frac{t^2}{8 c^2 k (\mathbb{E}(X)+t)}\right).
\]
\end{theorem}

\medskip 
We additionally recall two classical results concerning matchings.  
We include them for clarity of notation.

\begin{theorem}[Hall's marriage theorem~\cite{Hall}]\label{HallTh}
Let $G=(V,E)$ be bipartite with bipartition $V=X\cup Y$.  
Then $G$ contains a matching saturating $X$ if and only if
\[
\forall S\subseteq X:\quad |N(S)|\ge |S|.
\]
\end{theorem}

\begin{theorem}[Berge's theorem~\cite{Berge}]\label{BergeTh}
A matching $M$ in a graph $G$ is maximum if and only if $G$ contains no augmenting path with respect to~$M$.
\end{theorem}

We shall finally also make use of the following 
tool, originating from the theory of majority edge-colourings~\cite{PP},
which bounds a certain discrepancy parameter.  
Given a $k$-edge-colouring $\omega:E\to[k]$ of a graph $G=(V,E)$, set $E_i=\omega^{-1}(i)$ for $i\in[k]:=\{1,2,\ldots,k\}$. Let us denote
\[
\mathcal{D}_k(G):=\min_{\omega}\; \max_{v\in V(G)}\max_{1\le i<j\le k}
       \big|\, d_{E_i}(v)-d_{E_j}(v)\,\big|,
\]
where $\omega$ ranges over all $k$-edge-colourings of $G$.

\begin{theorem}[Przyby{\l}o--P\k{e}ka{\l}a, Observation~29 in~\cite{PP}]\label{PP}
For every graph $G$ and every integer $k\ge 2$, 
$
\mathcal{D}_k(G)\le 2.
$
\end{theorem}
As a direct consequence, we obtain the following corollary
\begin{corollary}\label{PP-col}
Every graph $G$ can be edge-decomposed into $k$ spanning subgraphs 
$G_1,\dots,G_k$ such that for every vertex $v$ and every $i\in[k]$,
\[
\frac{d_G(v)}{k}-2 \le d_{G_i}(v) \le \frac{d_G(v)}{k}+2.
\]
\end{corollary}

This balanced decomposition shall play a key role in controlling degrees throughout our main lemma below.

\section{Main Lemma}\label{MainLemmaSection}

The main purpose of the lemma below is to satisfy the vast majority of edges of $G$. 
To that end we  decompose $G$ in several steps into $G'$, containing almost all edges of $G$, and a leftover graph $G''$, to be used later on to complete the proofs of the main results of the paper.
In order to satisfy most of the edges we shall colour only a tiny part of the edges in $G'$, thus constructing a (proper) partial edge-colouring of $G$.
To achieve this, we shall first specify a series $G_1,\ldots,G_s$ of edge-disjoint subgraphs of $G'$, and then we shall randomly select a matching to be coloured with a single colour in each of these subgraphs.
The obtained colouring shall satisfy our requirements.
A limited number of still unsatisfied edges shall be then adjoined to the edges of $H''$ -- a subgraph extracted from $G$ at the beginning of our argument, to form $G''$. The latter graph shall encapsulate all unsatisfied edges and large enough number of yet not coloured edges (from $H''$) to admit 
completing proofs of the main results of the paper,
both in non-proper and proper settings.

\begin{lemma}\label{ML}
There exists $\Delta_0$ such that for every graph $G$ with maximum degree $\Delta \geq \Delta_0$ and minimum degree 
$ 
\delta \geq \Delta - 2\sqrt{\Delta}\,\log_2 \Delta,
$ 
there is a partial edge-colouring and an edge-decomposition of $G$ into $G'$ and $G''$ in which:
\begin{enumerate}
  \item each colour class forms a matching,
  \item at most $\log_2 \Delta + 1$ colours are used,
  \item all edges of $G'$ are satisfied,
  \item for every $v\in V$, $d_{G''}(v) \leq 4\ln\Delta$,
  \item $G''$ contains a set $U$ of uncoloured edges such that for every $v\in V$,
    \begin{equation}\label{dU2ineq}
    289\ln\ln\Delta \leq d_{U}(v) \leq 290\ln\ln\Delta.
  \end{equation}
\end{enumerate}
\end{lemma}

\begin{pf}
Let $G=(V,E)$ be a graph with maximum degree $\Delta$ and minimum degree 
\begin{equation}\label{deltalarge}
\delta \geq \Delta - 2\sqrt{\Delta}\,\log_2 \Delta.
\end{equation}
We also assume that $\Delta$ is large enough so that all explicit inequalities below hold.

\subsection{Graph Decompositon}

We begin by applying Corollary~\ref{PP-col} to the graph $G$ with 
$k = \lceil\frac{\Delta}{290\ln\ln\Delta-2}\rceil$. 
As a result, we obtain subgraphs $H_1, \ldots, H_k$ of $G$. 
Let us set $H'' = H_1$ and $H' = G - H_1=(V,E\smallsetminus E(H_1))$. 
By Corollary~\ref{PP-col} and~(\ref{deltalarge})
it follows that for every $v\in V$, 
\begin{align} 
 d_{H''}(v) \; & \leq\; \frac{\Delta}{k}+2 \; \leq \; 290 \ln\ln \Delta, \label{dH''est1}\\
d_{H''}(v)  \; & \geq\; \frac{\delta}{k}-2 \; \geq \; \frac{\Delta - 2\sqrt{\Delta}\,\log_2 \Delta}{k}-2 \; \geq \; 289 \ln\ln \Delta \label{dH''est2}
\end{align}
for $\Delta$ large enough, and hence,
\begin{equation}\label{dH'estx2}
d_G(v) - 290\ln\ln \Delta \;\leq\; d_{H'}(v) \;\leq\; d_G(v) - 289\ln\ln \Delta.
\end{equation}

As mentioned before, we shall set aside $H''$ and use later after providing 
a partial colouring for the remaining graph.  
Before doing so, let us apply Corollary~\ref{PP-col} 
once more -- this time to the graph $H'$,
thereby decomposing it into 
\[
s = \lceil \log_2 \Delta \rceil
\]
edge-disjoint subgraphs $G_1, \ldots, G_s$. 
By  Corollary~\ref{PP-col}
and~(\ref{dH'estx2}), 
for every vertex $v\in V$ 
and each $i\in[s]$ we have 
\begin{align}
 d_{G_i}(v) &\ge \frac{d_{H'}(v)}{s}-2 
\ge\frac{\Delta-2\sqrt{\Delta}\log_2\Delta-290\ln\ln\Delta}{s}-2\ge\frac{\Delta}{s}-3\sqrt{\Delta},
\nonumber\\
d_{G_i}(v) &\leq \frac{d_{H'}(v)}{s}+2
\le\frac{\Delta-289\ln\ln\Delta}{s}+2\le\frac{\Delta}{s}+2, \nonumber
\end{align}
which implies that for each $i\in [s]$, 
\begin{equation}\label{ReqDeltaBounds}
\delta_i:=\delta(G_i)\;\ge\;\frac{\Delta}{s}-3\sqrt{\Delta}
\qquad\text{and}\qquad
\Delta_i:=\Delta(G_i)\;\le\;\frac{\Delta}{s}+2. 
\end{equation}

\subsection{Vertex partitions}
In every graph $G_i$ we shall now extract an appropriate matching in a randomized manner. 
Instead of directly choosing the matching, for each $i\in[s]$ we first generate a random 
partition of the vertex set into two almost equal parts with certain additional features. 
These properties ensure, due to Hall’s and Berge’s theorems, that within the first part 
one can find a large matching, covering nearly all its vertices. For technical convenience this 
first part shall be split further into three subsets, denoted $V_1^{(i)},V_2^{(i)},V_3^{(i)}$, where $V_3^{(i)}$ is supposed to be very small. 
The remaining vertices shall form $V_4^{(i)}$. 
A matching $M_i$ chosen inside $V_1^{(i)} \cup V_2^{(i)} \cup V_3^{(i)}$, and coloured with $i$, shall saturate all vertices of $V_1^{(i)} \cup V_2^{(i)}$. As a consequence, every edge between $V_1^{(i)} \cup V_2^{(i)}$ and $V_4^{(i)}$ shall have exactly one endpoint incident with a colour-$i$ edge, and shall therefore be satisfied.
By randomness, with high probability each vertex shall get about half 
of its incident edges satisfied in this way. Since our random choices for different $i\in[s]$ 
shall be independent, the collection of $s$ such matchings shall be expected to satisfy in total the vast majority 
of all edges of $H'$.

For every $i\in [s]$ we set
\begin{equation}\label{EpsDef}
\varepsilon_i=\sqrt{\Delta_i}\,\ln \Delta_i,
\end{equation}
noting it is much larger than $\sqrt{\Delta}$, for $\Delta$ large enough.  
For every vertex $v\in V$ we introduce independent random variables $R_{v,i}$, $i\in [s]$, defined by
\begin{equation}\label{Zdef}
R_{v,i}=\left\{\begin{array}{lll}
1 & {\rm with~probability} & 1/4 - 3\varepsilon_i/\Delta_i\\ 
2 & {\rm with~probability} & 1/4\\ 
3 & {\rm with~probability} & 6 \varepsilon_i/\Delta_i\\ 
4 & {\rm with~probability} & 1/2 - 3 \varepsilon_i/\Delta_i 
\end{array}\right..
\end{equation}
These random variables define the intended vertex partition for each $G_i$, namely
$$
V_q^{(i)}=\{v\in V: R_{v,i}=q\}, \qquad q\in[4].
$$
Thus every $G_i$ comes with a partition $V=V_1^{(i)}\cup V_2^{(i)}\cup V_3^{(i)}\cup V_4^{(i)}$.  

For a vertex $v\in V$ we additionally set
\[
R_{v,1}^{(i)}=|N_{G_i}(v)\cap V_1^{(i)}|,\qquad
R_{v,2}^{(i)}=|N_{G_i}(v)\cap V_2^{(i)}|,\qquad
R_{v,13}^{(i)}=|N_{G_i}(v)\cap (V_1^{(i)}\cup V_3^{(i)})|.
\]
We shall be interested in the following events for each $v\in V$ and $i\in [s]$:
\begin{itemize}
    \item $D_{v,1}^{(i)} :~ R_{v,1}^{(i)} < \frac{\Delta_i}{4} - 2\varepsilon_i$,
    \item $F_{v,2}^{(i)} :~ R_{v,2}^{(i)} > \frac{\Delta_i}{4} - 2\varepsilon_i$,
    \item $F_{v,13}^{(i)} :~ R_{v,13}^{(i)} > \frac{\Delta_i}{4} + \varepsilon_i$,
    \item $D_{v,2}^{(i)} :~ R_{v,2}^{(i)} < \frac{\Delta_i}{4} + \varepsilon_i$.
\end{itemize}
In what follows, we shall show that the probability of the failure of any of these events is miniscule.

Note first that by~(\ref{Zdef}) we have
\(\mathbf{E}(R_{v,1}^{(i)}) \le \frac{\Delta_i}{4} - 3\varepsilon_i\).
Applying the Chernoff Bound together with~(\ref{EpsDef}) and~(\ref{ReqDeltaBounds}), we obtain
\begin{eqnarray}
    \mathbf{P}\left(\overline{D_{v,1}^{(i)}}\right) &=& \mathbf{P}\Big(R_{v,1}^{(i)} \ge \Big(\frac{\Delta_i}{4}-3\varepsilon_i\Big)+\varepsilon_i\Big)\nonumber\\
&\le& \exp\!\left(-\frac{\varepsilon_i^2}{3\left(\frac{\Delta_i}{4}-3\varepsilon_i\right)}\right)
\le \exp\!\left(-\frac{4}{3}(\ln\Delta_i)^2\right)
\le \Delta^{-3}. \label{Sv1Bound}
\end{eqnarray}

Similarly, by~(\ref{Zdef}), (\ref{EpsDef}) and~(\ref{ReqDeltaBounds}) we have
\(\mathbf{E}(R_{v,2}^{(i)}) \ge \frac{\delta_i}{4} \ge \frac{\Delta_i}{4} - \varepsilon_i\), and
hence,
\begin{eqnarray}
\mathbf{P}\left(\overline{F_{v,2}^{(i)}}\right)
&=& \mathbf{P}\Big(R_{v,2}^{(i)} \le \Big(\frac{\Delta_i}{4}-\varepsilon_i\Big)-\varepsilon_i\Big)\nonumber\\
&\leq& \exp\!\left(-\frac{\varepsilon_i^2}{2\left(\frac{\Delta_i}{4}-\varepsilon_i\right)}\right)
\le \exp\!\left(-2(\ln\Delta_i)^2\right)
\le \Delta^{-3}. \label{Bv2Bound}
\end{eqnarray}

Further, again by~(\ref{Zdef}), (\ref{EpsDef}) and~(\ref{ReqDeltaBounds}),
$ 
\mathbf{E}(R_{v,13}^{(i)}) 
\ge \delta_i\Big(\frac{1}{4}+\frac{3\varepsilon_i}{\Delta_i}\Big)
\ge \frac{\Delta_i}{4}+2\varepsilon_i,
$ 
and hence,
\begin{eqnarray}
\mathbf{P}\left(\overline{F_{v,13}^{(i)}}\right)
&=& \mathbf{P}\Big(R_{v,13}^{(i)} \le \Big(\frac{\Delta_i}{4}+2\varepsilon_i\Big)-\varepsilon_i\Big)\nonumber\\
&\le& \exp\!\left(-\frac{\varepsilon_i^2}{2\left(\frac{\Delta_i}{4}+2\varepsilon_i\right)}\right)
\le \exp\!\left(-(\ln\Delta_i)^2\right)
\le \Delta^{-3}. \label{Bv13Bound}
\end{eqnarray}

Finally, by~(\ref{Zdef}) we have \(\mathbf{E}(R_{v,2}^{(i)}) \le \frac{\Delta_i}{4}\), hence
\begin{eqnarray}
\mathbf{P}\left(\overline{D_{v,2}^{(i)}}\right)
&=& \mathbf{P}\Big(R_{v,2}^{(i)} \ge \frac{\Delta_i}{4}+\varepsilon_i\Big)\nonumber\\
&\le &\exp\!\left(-\frac{\varepsilon_i^2}{3\cdot\frac{\Delta_i}{4}}\right)
= \exp\!\left(-\frac{4}{3}(\ln\Delta_i)^2\right)
\le \Delta^{-3}. \label{Sv2Bound}
\end{eqnarray}

We shall later explain how the events
\(D_{v,1}^{(i)},\, F_{v,2}^{(i)},\, F_{v,13}^{(i)},\, D_{v,2}^{(i)}\)
collectively guarantee the existence of the required large matchings.

For each $i\in[s]$, every vertex $v$ independently chooses one of
the four sets $V_1^{(i)},\dots,V_4^{(i)}$ according to the distribution
described earlier.  
For a fixed vertex $v$, we encode the entire sequence of its $s$
choices by a vector
\[
w = (w_1,\dots,w_s)\in[4]^s,
\]
where $w_i=x$ means $v\in V_x^{(i)}$, and call it a \emph{vector of draws}.

We say an edge $vu\in E$ is in \emph{good configuration} in partition $i$ if
\[
v\in V_{12}^{(i)}:=V_1^{(i)}\cup V_2^{(i)}
\quad\text{and}\quad
u\in V_4^{(i)},
\]
or symmetrically $u\in V_{12}^{(i)}$ and $v\in V_4^{(i)}$.
Otherwise, $vu$ is in \emph{bad confirguration} in the corresponding partition.
Finally, we say that $vu$ is \emph{always bad} if it lands in bad configuration in all $s$ partitions.
Let
\[
B_v=\bigl|\{u\in N_G(v): vu \text{ is always bad}\}\bigr|
\]
be the number of neighbours of $v$ in $G$ whose incident edges joining them with $v$ 
are in bad configurations in all partitions.
Our goal shall now be bounding the probability of the complement of the event:
\[
A_v: ~ B_v\leq 3\ln\Delta -3.
\]

The sets $V_3^{(i)}$ serve a technical role in our randomized construction, but are rather annoying  
necessities. These are however likely not to include many vertices. In order to control their influence, complicating further calculations, we shall explicitly separate some inconvenient,  
yet highly improbable vectors of draws.
For a given $w\in[4]^s$, let
\[
|w^{-1}(3)|=\bigl|\{i\in[s]: w_i=3\}\bigr|
\]
be the number of partitions $i$ in which a vertex $v$ with draw vector $w$ falls into $V_3^{(i)}$.
Set
\[
W:=\{w\in[4]^s: |w^{-1}(3)|\le 6\}.
\]
For each $w\in[4]^s$ and vertex $v\in V$, let $C_{v,w}$ be the event that the draws for
$v$ are precisely given by $w$, i.e., 
\[
C_{v,w}: ~w=(w_1,\dots,w_s) \text{ is the vector of draws of }v.
\]
Now define
\[
C_v=\bigcup_{w\notin W} C_{v,w}
\]
to be the event that $v$ falls into $V_3^{(i)}$ in at least $7$ partitions. 
Using the law of total probability we split \(\Pr(\overline{A_v})\) as follows:
\begin{equation}\label{TotalSplit}
\Pr(\overline{A_v})
= \Pr(\overline{A_v}\mid C_v)\Pr(C_v) + \Pr(\overline{A_v}\mid \overline{C_v})\Pr(\overline{C_v})
\le \Pr(C_v) + \Pr(\overline{A_v}\mid \overline{C_v}).
\end{equation}

We first bound $\Pr(\overline{A_v}\mid \overline{C_v})$. Fix any $w\in W$.  
Conditioning on $C_{v,w}$, the vertex $v$ belongs to $V_3^{(i)}$ in at most six partitions. 
Hence in at least $s-6$ partitions 
it lands in $V_1^{(i)}\cup V_2^{(i)}\cup V_4^{(i)}$, 
and for every corresponding $i$, we have: 
\[
\Pr\left(\text{$vu$ is in bad configuration in partition 
$i$}\mid C_{v,w}\right)
= \frac12 + \frac{3\varepsilon_i}{\Delta_i}
= \frac12 +\frac{3\ln\Delta_i}{\sqrt{\Delta_i}}
\le \frac12 + \frac{\ln^2\Delta}{\sqrt{\Delta}},
\]
for sufficiently large \(\Delta\). Thus, the probability that $vu$ is in bad configurations in all the relevant partitions 
(thus also the probability that $vu$ is always bad) 
is at most
\begin{equation}\label{SingleBad}
\left(\tfrac12 + \frac{\ln^2\Delta}{\sqrt{\Delta}}\right)^{s-6}.
\end{equation}
Note that such events are independent for all $u\in N_G(v)$ under the given conditioning. Moreover, 
since $s=\lceil\log_2\Delta\rceil$ and  $1+x\le e^x$ for every real number $x$, by~(\ref{SingleBad}),
\begin{align}
\mathbb{E}(B_v\mid C_{v,w})
&\le
\Delta\Bigl(\tfrac12 + \frac{\ln^2\Delta}{\sqrt{\Delta}}\Bigr)^{s-6}
= 
\Delta\Bigl(\tfrac12\Bigr)^{s-6} 
\Bigl(1+2\frac{\ln^2\Delta}{\sqrt{\Delta}}\Bigr)^{s-6}\nonumber\\
&\le \Delta \Bigl(\tfrac12\Bigr)^{\log_2\Delta-6} 
\left(e^{2\frac{\ln^2\Delta}{\sqrt{\Delta}}}\right)^{s-6}
= 64 e^{\frac{2(s-6)\ln^2\Delta}{\sqrt{\Delta}}} \leq 65
\end{align}
for large enough $\Delta$. 
By the Chernoff bound (Lemma~\ref{Ch2}), with
$t:=3\ln\Delta-68$, we thus obtain: 
\begin{align}\label{Chernoff_U}
\Pr\left(\overline{A_v}\mid C_{v,w}\right) =
\Pr\left(B_v>3\ln\Delta -3\mid C_{v,w}\right)
&\le \Pr\left(B_v\ge \mathbb{E}(B_v\mid C_{v,w})+t \mid C_{v,w}\right)\nonumber\\
&\le \exp\!\left(-\frac{(3\ln\Delta-68)^2}{3\ln\Delta-68+130}\right)\nonumber\\
&\le \exp\!\big(-2.5\ln\Delta\big)
= \Delta^{-2.5}, \nonumber
\end{align}
for all sufficiently large \(\Delta\). Using the law of total probability we thus get
\begin{eqnarray}
\Pr(\overline{A_v}\mid \overline{C_v})
&=&\sum_{w\in W} \Pr(\overline{A_v}\mid C_{v,w})\Pr(C_{v,w}\mid \overline{C_v})\nonumber\\
&\le& \Delta^{-2.5}\sum_{w\in W}\Pr(C_{v,w}\mid \overline{C_v})
~=~ \Delta^{-2.5}.\label{Av1}
\end{eqnarray}

It remains to bound \(\Pr(C_v)\). Let \(S_v\) be the number of
partitions 
in which \(v\in V_3^{(i)}\). Since for every partition, i.e. for every fixed $i\in[s]$, 
\[
\Pr\left(v\in V_3^{(i)}\right) = \frac{6\varepsilon_i}{\Delta_i}
\le \frac{\ln^2\Delta}{\sqrt{\Delta}}
\]
(for all sufficiently large \(\Delta\)), we obtain by a simple binomial bound:
\begin{align}\label{BvCalc}
\Pr(C_v)=\Pr(S_v\ge7)
&\le \binom{s}{7} \left(\frac{\ln^2\Delta}{\sqrt{\Delta}}\right)^7
\le s^7\left(\frac{\ln^2\Delta}{\sqrt{\Delta}}\right)^{7}\nonumber\\
&\le (\log_2\Delta+1)^7\cdot (\ln\Delta)^{14}\cdot \Delta^{-3.5} \leq \Delta^{-3}.
\end{align}
By~(\ref{TotalSplit}), (\ref{Av1}) and~(\ref{BvCalc}) we thus obtain that
\begin{equation}\label{Av-Bound}
\Pr(\overline{A_v})\le \Delta^{-3}+\Delta^{-2.5}\le \Delta^{-2.4}.
\end{equation}

Finally, observe that each event 
$\overline{D_{v,1}^{(i)}}$, $\overline{F_{v,2}^{(i)}}$, $\overline{F_{v,13}^{(i)}}$,
    $\overline{D_{v,2}^{(i)}}$ and $\overline{A_v}$
    depends only on the random
choices made within distance at most $1$ of \(v\),
hence it is mutually independent of all other such events associated with vertices $u$ which are at distance greater than $2$ from $v$.
Thus, each of these events is mutually independent of all but at most $\Delta^2(4s+1)$ other such events.
Moreover, by (\ref{Sv1Bound}), (\ref{Bv2Bound}), (\ref{Bv13Bound}), (\ref{Sv2Bound})  and~(\ref{Av-Bound}),
the probability of each of these events is bounded above by $\Delta^{-2.4}$.
Since $e\Delta^{-2.4}(\Delta^2(4s+1)+1) < 1$ for large enough $\Delta$, by the Lov\'asz Local Lemma, 
all events $D_{v,1}^{(i)}$, $F_{v,2}^{(i)}$, $F_{v,13}^{(i)}$,
    $D_{v,2}^{(i)}$ and $A_v$ hold simultaneously with positive probability
    -- let us fix any family of 
  partitions  $V=V_1^{(i)}\cup V_2^{(i)}\cup V_3^{(i)}\cup V_4^{(i)}$, $i\in[s]$,
  witnessing this.

\subsection{Existence of suitable matchings}

We now argue that in each fixed $G_i$ with $i\in [s]$ one can select a matching $M_i$ whose
edges lie entirely in $V_1^{(i)}\cup V_2^{(i)}\cup V_3^{(i)}$ and which
covers all vertices of $V_1^{(i)}\cup V_2^{(i)}$.  
Consider first the bipartite graph
\[
B_i = G_i[V_1^{(i)},V_2^{(i)}].
\]
By the events $F_{u,2}^{(i)}$ and $D_{v,1}^{(i)}$, every
$u\in V_1^{(i)}$ has larger degree in $B_i$ than any $v\in V_2^{(i)}$, so a standard
double-counting argument implies Hall’s condition.  
Consequently, Theorem~\ref{HallTh} yields a matching
$M_i^{\mathrm{base}}$ in $B_i$ which saturates
$V_1^{(i)}$.

Now consider the extension of $B_i$ to the bipartite graph
$B'_i = G_i[V_1^{(i)}\cup V_3^{(i)},\, V_2^{(i)}].$
By the events  $F_{u,13}^{(i)}$ and $D_{v,2}^{(i)}$ we have
\[
d_{B'_i}(u) > d_{B'_i}(v)
\qquad\text{for all } \; \; u\in V_2^{(i)},\; v\in V_1^{(i)}\cup V_3^{(i)}.
\]
Hence, Hall’s condition also holds for $B'_i$, and so $B'_i$ contains a matching of size $|V_2^{(i)}|$.
Note $M_i^{\mathrm{base}}$ is also a matching in $B'_i$.
If it is not maximum in $B'_i$, then by Berge’s theorem
there exists an augmenting path with respect to $M_i^{\mathrm{base}}$.  
Flipping the matching adges along such a path produces a strictly larger matching that still
saturates all vertices of $V_1^{(i)}$.
Iterating this augmentation process eventually yields a maximum matching $M_i$
in $B'_i$, which by construction saturates both $V_1^{(i)}$ and $V_2^{(i)}$.  
Such $M_i$ is a desired matching.
For every $i\in[s]$ we fix such a matching $M_i$ and colour all its edges
with colour $i$.  
Let $c$ be the resulting partial edge-colouring of $G$.

Consider any vertex $v\in V$.  
By the event $A_v$, at most $3\ln\Delta-3$ edges incident with $v$ are
\emph{always bad}, that is, in bad configuration in all fixed $s$ partitions of $V$.
Hence, at least $d_G(v)-(3\ln\Delta-3)$
edges incident with $v$ are (each) 
in good configuration in at least one partition.
Consider any such edge $uv$ and the corresponding $i$-th partition.
Thus, one endpoint of $uv$
belongs to $V_{12}^{(i)}=V_1^{(i)}\cup V_2^{(i)}$ while the other belongs to
$V_4^{(i)}$. By construction the matching $M_i$ is chosen inside
$V_1^{(i)}\cup V_2^{(i)}\cup V_3^{(i)}$ and saturates all vertices of
$V_{12}^{(i)}$, whereas no vertex of $V_4^{(i)}$ is incident with an edge
of $M_i$. Consequently, exactly one endpoint of $uv$ is incident with a
colour-$i$ edge, and that colour therefore appears exactly once in the
 neighbourhood of $uv$. Hence, $uv$ is satisfied by colour $i$. 

Let $G'$ be the subgraph of $G$ induced by all the satisfied edges in $H'$,
and let $H'''=H'-G'$ be the subgraph of $H'$ induced by the remaining edges.
Let $G''=H''\cup H'''$ be the graph formed of the edges in $H''$, neither of which is coloured, 
and the edges in $H'''$, which are not satisfied. 
Clearly $G'$ and $G''$ form an edge-decomposition of $G$,
which is partially coloured with  $s\leq \log_2\Delta+1$ colours, each of which forms a matching, and all edges of $G'$ are satisfied.

Finally, in order to legitimise the last two requirement of the lemma,  
consider any $v\in V$. 
Let us denote $U=E(H'')$.
Note that  $d_{G''}(v)=d_{H''}(v)+d_{H'''}(v) = d_{U}(v)+d_{H'''}(v)$. By our construction, all edges in $U$ are uncoloured, and by~(\ref{dH''est1}) and (\ref{dH''est2}) we have, $289\ln\ln\Delta \leq d_U(v) \leq 290\ln\ln\Delta$.
Moreover, 
since $H'''$ includes unsatisfied edges exclusively, we have $d_{H'''}(v)\leq 3\Delta - 3$, which concludes the proof
for large enough $\Delta$.
\qed

\end{pf}

\section{Main result}\label{MainResultSection}

The following theorem shall immediately imply our two main results.
In order to avoid repetitions we impose within it a requirement which is not needed for one of these results.
Namely, we shall provide a \emph{proper} partial colouring of $G$, i.e. an assignment of colours to selected edges of $G$ such that no two adjacent edges are coloured the same.

\begin{theorem}\label{Main0}
There exists $\Delta_0$ such that 
every graph $G$ with maximum degree $\Delta \geq \Delta_0$ and minimum degree 
$\delta \geq \Delta - 2\sqrt{\Delta}\log_2\Delta$
admits a proper partial colouring $c$ with at most
\[
\log_2\Delta + 4000\ln\ln\Delta-1
\]
colours which satisfies all the edges of $G$.
\end{theorem}

\begin{pf}
We begin by applying Lemma~\ref{ML} to $G$.
As a result we obtain a partial colouring $c$ of $G$ with at most $\log_2\Delta+1$ colours and an edge-decomposition of $G$ to $G'$ and $G''$ such that  
all edges of $G'$ are already satisfied. 
Since the colour classes of $c$ form matchings, $c$ is clearly a proper partial colouring.
Consequently, we are left with the subgraph $G''$ of maximum degree $\Delta''=\Delta(G'')\leq 4\ln\Delta$
in which some edges may already be coloured, but none of them are satisfied. 
Our main goal shall be thus extending $c$ to a proper partial colouring of edges in $G''$, using new colours, so that all the edges in $G''$ get satisfied.
We shall be colouring only edges in the subset $U$ of the edges of $G''$, resulting from Lemma~\ref{ML}, none of which has been coloured yet.
We shall also abuse the adopted notation slightly, and we shall denote the colouring  of $G''$ to be constructed below by $c$ as well (treating it as an extension of the colouring of $G'$ to a disjoint set of edges and with a new set of colours). Let us recall that by~(\ref{dU2ineq}), for every $v\in V$,
\begin{equation}\label{RepeateddU-estx2}
289\ln\ln\Delta\leq d_U(v)\leq 290\ln\ln\Delta.
\end{equation}

Now, we colour the edges of $U$ randomly using
$M := \lfloor 3999\ln\ln\Delta\rfloor$ new
colours in the following manner: first each $e\in U$ receives a colour $c'(e)$ chosen uniformly and independently at random from the $M$ colours; next every two adjacent edges from $U$ coloured the same get both uncoloured.
The final resulting partial colouring of $U$ we denote by $c$.
Note that by the adopted uncolouring rule, $c$ is clearly proper.

For every fixed edge $e\in E(G'')$, we define the \emph{open neighbourhood of $e$ in $U$} by
\[
U_e := \{uv \in U : (u\in e \lor v\in e)\,\}\smallsetminus \{e\}; 
\qquad u_e:=|U_e|.
\]
Note that for any $e=vv'\in E(G'')$, since $(d_U(v)-1)+(d_U(v')-1)\leq d_U(e)\leq d_U(v)+d_U(v')$, then by~(\ref{RepeateddU-estx2}),
\begin{equation}\label{ue-estimate}
 576\ln\ln\Delta\leq u_e\leq  580 \ln\ln\Delta.
\end{equation}
For any $e\in E(G'')$ and $r\in U_e$, let us further define the indicator variable
\[
Y_{e,r} := \mathbf{1}\{\exists\,r'\in U_e\cup U_r \setminus\{r\} : c'(r')=c'(r)\}.
\]
Note that if $Y_{e,r} =0$ for at least one $r\in U_e$, then $e$ is satisfied in the resulting colouring $c$.
By~(\ref{RepeateddU-estx2}), since $e$ and $r$ have a common vertex,  
\begin{equation}\label{UeUr-bound}
|U_e\cup U_r \setminus\{r\}|\leq 3\cdot 290\ln\ln\Delta = 870\ln\ln\Delta.
\end{equation}
Set $Y_e = \sum_{r\in U_e} Y_{e,r}$.
Observe that if, for every $e \in E(G'')$, it holds that $Y_e < u_e$, 
then every edge $e \in E(G'')$ is satisfied in the final colouring $c$, as desired. For each $e\in E(G'')$ we shall thus bound the probability of the opposite event:
$$Q_e: Y_e=u_e.$$

Note first that for fixed \(e\) and \(r\in U_e\), by~(\ref{UeUr-bound}) we have
\[
\Pr(Y_{e,r}=1)\leq 1-\Big(1-\frac{1}{M}\Big)^{870\ln\ln\Delta}.
\]
Since $1-x\geq e^{-\frac{x}{1-x}}$ for $0< x<1$, we also have
\[
\Big(1-\frac{1}{M}\Big)^{870\ln\ln\Delta} \ge \exp\!\Big(-\frac{870\ln\ln\Delta}{M-1}\Big).
\]
Combining the two inequalities above, we thus obtain
\[
\Pr(Y_{e,r}=1)\le 1-\exp\!\Big(-\frac{870\ln\ln\Delta}{M-1}\Big)
\le 1-\exp\!\Big(-\frac{870}{3998}\Big)
< \frac{1}{5}.
\]
Therefore, by the linearity of expectation,
\[
\mathbb{E}(Y_e)=\sum_{r\in U_e}\mathbb{E}(Y_{e,r}) \le \frac{u_e}{5}.
\]

Let $T_e$ be the trial corresponding to the choice of colour $c'(e)$ for an edge $e\in U$. Thus, $Y_e$ is determined by $\{T_{r'} : r' \in \bigcup_{r\in U_e}U_r\}$. Note that changing the value of any $T_r$ may change $Y_e$ by at most $2$, and the fact that $Y_e \geq s$ can be certified by assignments of colours by $c'$ to (no more than) $2s$ edges. Thus, the assumptions of Theorem~\ref{tlg} are satisfied with 
$c=k=2$, and hence we apply it  below with $t = \frac{3u_e}{5}$.
Since $\mathbb{E}(Y_e)\leq \frac{u_e}{5}$, by~(\ref{ue-estimate}) we thus obtain 
\begin{align}
\Pr\left(Q_e\right) &=
\Pr\left(Y_e=u_e\right) \leq
\Pr\left(\left|Y_e - \mathbb{E}(Y_e)\right| > \frac{3u_e}{5} + 40\sqrt{2\mathbb{E}(Y_e)} + 512\right) \nonumber\\ 
&\le 4e^{- \frac{1}{64}\frac{\frac{9}{25}u_e^2}{(\mathbb{E}(Y_e) + \frac{3}{5}u_e)}}
\le 4e^{- \frac{1}{64}\frac{\frac{9}{25}u_e^2}{\frac{4}{5}u_e}}
= 4e^{- \frac{9u_e}{20\cdot 64}}
\leq 4e^{- \frac{9\cdot 576\ln\ln\Delta}{20\cdot 64}}
= 4\left(\ln\Delta\right)^{-4.05}.
\end{align}

Observe that each event $Q_e$ is mutually independent of the events $Q_{e'}$ with $e'$ at distance 
at least $5$ (in $L(G'')$), that is of all but at most $2(\Delta'')^4\leq 2(4\ln\Delta)^4$ of the considered events. 
For $\Delta$ large enough, by the Lov\'asz Local Lemma, we may thus infer that 
with positive probability, none of $Q_e$ holds, which means that there must exist a desired colouring $c$ such that $Y_e < u_e$ for every $e\in E(G'')$. 
By construction, taking into account the colours from $G'$, the obtained proper partial colouring of $G$ exploits no more than $\log_2\Delta+1+M\leq \log_2\Delta+4000\ln\ln\Delta-1$ colours, for $\Delta$ large enough.
\qed
\end{pf}

\section{Final results}\label{FinalResultsSection}

We now formally state our target theorems, which are direct corollaries of  Theorem~\ref{Main0}.

\begin{theorem}\label{Main1}
There exists $\Delta_0$ such that  
\[
\chi_{\rm OCF}'(G) \leq \log_2\Delta + 4000\ln\ln\Delta
\]
for every graph $G$ with maximum degree $\Delta \geq \Delta_0$ and minimum degree 
\(\delta \geq \Delta - 2\sqrt{\Delta}\log_2\Delta\).
\end{theorem}
\begin{pf}
By Theorem~\ref{Main0} there exists a partial edge colouring of $G$ with at most $\log_2\Delta + 4000\ln\ln\Delta-1$ colours which satisfies all its edges. To complete this it is enough to use a single additional colour and assign it to all
the remaining uncoloured edges. 
\qed
\end{pf}

\begin{theorem}\label{Main2}
There exists $\Delta_0$ such that  
\[
\chi_{\rm pOCF}'(G) \leq \Delta + 2\log_2\Delta
\]
for every graph $G$ with maximum degree $\Delta \geq \Delta_0$ and minimum degree 
\(\delta \geq \Delta - 2\sqrt{\Delta}\log_2\Delta\).
\end{theorem}

\begin{pf}
As above,
by Theorem~\ref{Main0} there exists a partial \emph{proper} edge colouring of $G$ with at most $\log_2\Delta + 4000\ln\ln\Delta-1$ colours which satisfies all its edges. 
Clearly we may complement it to an open conflict-free proper edge-colouring of $G$ exploiting to that end at most $\Delta+1$ new colours, thus using in total no more than $\Delta + \log_2\Delta + 4000\ln\ln\Delta \le \Delta + 2\log_2\Delta$ colours, for $\Delta$ large enough.
\qed
\end{pf}

By Theorems~\ref{Main1} and \ref{Main2}, due to~(\ref{Eq1.2}) and~(\ref{Eq1.3}), we thus obtain the following conclusion.

\begin{corollary}\label{MainCor12}
for every graph $G$ with maximum degree $\Delta$ and minimum degree 
$\delta \geq \Delta - 2\sqrt{\Delta}\log_2\Delta$,
\begin{align}
\chi_{\rm OCF}'(G) &\leq \log_2\Delta + O(\ln\ln\Delta) = (1+o(1))\log_2\Delta,\nonumber\\
\chi_{\rm pOCF}'(G) &\leq \Delta + O(\ln\Delta) =  (1+o(1))\Delta. \nonumber
\end{align}
\end{corollary}

Assumptions in our results above were designed to admit their direct extensions to the random graphs setting.
In order to derive the corresponding theorems, we need the following technical observation on random graphs
from~\cite{KPreg}. 

\begin{observation}[\cite{KPreg}, Observation~11]\label{gnp_obs}
    If $G = G(n,p)$ is a random graph with $p \gg n^{-\varepsilon}$ for some fixed
    constant $\varepsilon \in (0,1)$, then asymptotically almost surely,
    \[
        \delta(G) \;\ge\; \Delta(G) - 2\sqrt{\Delta(G)}\big(\ln\Delta(G)\big)^{3/4}
        \qquad\text{and}\qquad
        \Delta(G)\gg n^{\,1-\varepsilon}.
    \]
\end{observation}

Since $\Delta - 2\sqrt{\Delta}(\ln\Delta)^{3/4}\geq \Delta - 2\sqrt{\Delta}\log_2\Delta$ for $\Delta$ large enough,
Corollary~\ref{MainCor12} combined with Observation~\ref{gnp_obs} immediately imply the following results.

\begin{theorem}\label{gnp-thm-open}
    If $G = G(n,p)$ is a random graph with $p \gg n^{-\varepsilon}$ for some fixed
    constant $\varepsilon \in (0,1)$, then asymptotically almost surely,
    \[
        \chi_{\rm OCF}'(G)
        \;\le\;
        \log_2\Delta(G) + O\!\big(\ln\ln\Delta(G)\big)
        \;=\;
        (1+o(1))\,\log_2\Delta(G).
    \]
\end{theorem}

\begin{theorem}\label{gnp-thm-proper}
If $G = G(n,p)$ is a random graph with $p \gg n^{-\varepsilon}$ for some fixed
    constant $\varepsilon \in (0,1)$, then asymptotically almost surely,
    \[
        \chi_{\rm pOCF}'(G)
        \;\le\;
        \Delta(G) + O\!\big(\ln\Delta(G)\big)
        \;=\;
        (1+o(1))\,\Delta(G).
    \]
\end{theorem}

\section{Concluding remarks}\label{ConcludingRemarksSection}

Multiplicative constants at smaller order terms in upper bounds we obtained, in particular in Theorems~\ref{Main1} and~\ref{Main2},  are surely not optimal, even within the framework of proving techniques we used.
One reason for this was the aimed clarity of presentation of our reasonings. 
Moreover, in order to construct a universal tool, formulated within 
Theorem~\ref{Main0}, applicable to both $\chi_{\rm OCF}'(G)$ and $\chi_{\rm pOCF}'(G)$,
we imposed properness of the colouring it provides.
Dropping this requirement would surely allow to reduce the constant $4000$ in Theorems~\ref{Main1} significantly.
Our main goal was however optimising leading constants and proving asymptotically optimal bounds, formulated in  Corollary~\ref{MainCor12}.

In order to see that the bound for $\chi_{\rm OCF}'(G)$ in Corollary~\ref{MainCor12} is indeed asymptotically optimal, in particular for regular graphs, it suffices to consider the complete graphs, $K_n$.
To show these require $\log_2\Delta - O(\ln\ln\Delta) = (1-o(1))\log_2\Delta$ colours, we apply a reasoning very similar to the one used in~\cite{DebskiPrzybylo} (Theorem~3).

\begin{observation}
For every $n\geq 3$, $\chi_{\rm OCF}'(K_n)> \log_2 n - \log_2\log_2 n -1$.
\end{observation}

\begin{pf}
We argue by contradiction. So suppose there exists an open conflict-free edge colouring $c$ of $K_n$ with 
$k:=\lfloor \log_2 n - \log_2\log_2 n -1\rfloor$ colours. We may assume that $n\geq 16$; otherwise $k\leq 0$. For every vertex $v$, let $P(v)$ denote the \emph{palette} of colours incident with $v$, i.e. the set of all colours assigned to the edges incident with $v$. 
Since no more than $2^k$ distinct palettes may occur, at least one of them, say $Q$, is associated with at least $\frac{n}{2^k}$ vertices -- denote the set of all these vertices by $S$. Set $s=|S|$. 
Thus, 
\begin{eqnarray}\label{s-bound}
s\geq \frac{n}{2^k}\geq 2\log_2 n > \frac{4}{3}k+1.
\end{eqnarray} 

Consider the set of edges $F$ induced by $S$, each of which is satisfied -- for every such edge $e$ fix one colour $\omega(e)\in Q$ which occurs exactly once in the open neighbourhood of $e$. 
For any given colour $a\in Q$, denote by $F_a$ the set of all edges $e\in F$ with $\omega(e)=a$.
Since for any such edge $e=uv\in F_a$, we have $P(u)=Q=P(v)$, the colour $a$ must be assigned to $e$ and to exactly one edge adjacent to $e$. Consequently, every connected component induced by the edges in any given $F_a$ must either be an isolated edge or a path with two edges, and thus $|F_a|\leq \frac{2}{3}s$. 
Therefore, the total number of edges in all the sets $F_a$ does not exceed
$$\frac{2s}{3}k = \frac{s}{2} \cdot \frac{4k}{3} <\frac{s}{2}(s-1),$$
where the inequality follows by~(\ref{s-bound}).
Since the sets $F_a$ partition $F$ and $|F|={s\choose 2}$, we obtain a contradiction.
\qed
\end{pf}

For general graphs, apart from the bound $\chi_{\rm OCF}'(G)=O(\ln^{2+\varepsilon})$, stemming from~(\ref{Eq1.2}), not much more is known. Small degree vertices are potentially problematic on the way towards adapting techniques developed in several previous related papers to the setting of open neighbourhoods. One of the issues encompasses potential isolated edges emerging while using approaches based on sequential graph decompositions. 
We believe the following question is one of the most natural and important problems to resolve in the open setting.
\begin{problem}
Is it true that  $\chi_{\rm OCF}'(G)=O(\log_2\Delta)$ for all graphs $G$ with maximum degree $\Delta$?
\end{problem}

As for the proper colouring setting, an obvious trivial lower bound $\chi_{\rm pOCF}'(G)\geq \Delta$ holds for every graph. Thus the corresponding result from Corollary~\ref{MainCor12} is asymptotically optimal.
Intuitively, proper colourings seem easier to handle towards obtaining an asymptotically optimal bound for general graphs. 
Note in particular that for a given properly edge-coloured graph $G$ and its edge $e=uv$,
if only $d_G(u)\neq d_G(v)$, 
then there must exist a unique colour in the open neighbourhood of $e$. 
So only adjacent vertices of the same degrees require extra care. 
Thus, in a way, graphs close to regular are the most challenging ones.
This is however precisely the regime where our Theorem~\ref{Main2} provides a meaningful estimate for a substantial family of graphs (in fact the majority of graphs, as the Observation~\ref{gnp_obs} on random graphs implies).
For general graphs, only an upper bound slightly above $2\Delta$ can be inferred from~(\ref{Eq1.3}), which seems to be very far from tight. In fact,
we are prone to believe that one should be able to improve it to $\Delta+C$, for some absolute constant $C$,
and
we pose the following conjecture.
\begin{conjecture}\label{Conjecture-proper}
For every graph $G$ with maximum degree $\Delta$,
$\chi_{\rm pOCF}'(G)\leq \Delta+3$. 
\end{conjecture} 
Note that if this conjecture holds, then the cycle $C_5$ proves it cannot be improved in general.


\begin{thebibliography}{99}
\bibitem{DATA}
Data sharing not applicable to this article as no datasets were generated or analysed during the current study
\end{thebibliography}
\end{document}